\documentclass{article}

\usepackage{epsfig}
\usepackage{graphicx}
\usepackage{amssymb}
\usepackage{psfrag}
\usepackage{amsfonts}
\usepackage{amsmath}
\usepackage{subfigure}
\usepackage[dvips]{color}

\newtheorem{lemma}{Lemma}

\newtheorem{theorem}[lemma]{Theorem}
\newtheorem{corollary}[lemma]{Corollary}
\newtheorem{conjecture}[lemma]{Conjecture}

\def\pf{\noindent {\bf Proof~~}}
\def\endpf{\hfill\rule{2mm}{2mm}\\\bigskip}

\title{Every longest circuit of a $3$-connected, $K_{3,3}$-minor free graph has a chord}
\author{Etienne Birmel\'e\\
Laboratoire Statistique et G\'enome, 
Universit\'e d'Evry,\\
UMR CNRS 8071, INRA 1152,\\
Tour Evry 2, 523 place des Terrasses de l'Agora,\\
91034 Evry Cedex,\\
France\\
{\tt etienne.birmele@genopole.cnrs.fr}}

\begin{document}

{\em This is a preprint of an article accepted for publication in Journal of Graph Theory \copyright  copyright 2008 
Wiley Periodicals }

\maketitle

\begin{abstract}
Carsten Thomassen conjectured that every longest circuit in a $3$-connected 
graph has a chord. We prove the conjecture for graphs having no $K_{3,3}$ minor, and consequently
for planar graphs.
\end{abstract}

Carsten Thomassen made the following conjecture~\cite{Bon95,Tho84}, where a circuit denotes a connected $2$-regular graph:

\begin{conjecture}[Thomassen]
Every longest circuit of a $3$-connected graph has a chord.
\end{conjecture}

That  conjecture has been proved for  planar 
graphs with minimum degree at least four~\cite{Zha87}, 
cubic graphs~\cite{Tho97}
 and graphs embeddable in several 
surfaces~\cite{KNZ07,LZ03,LZ03bis}. In this paper, we prove it for planar 
graphs in general. In fact, our result concerns a class of graphs 
which contains the 
planar graphs.
\medskip
 
\begin{figure} [htb]
\begin{center}
\epsfig{file=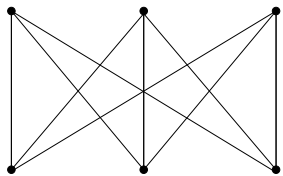, width=.3\textwidth}
\end{center}
\caption{The graph $K_{3,3}$}
\label{K33}
\end{figure}

Let us denote by $K_{3,3}$ the complete bipartite graph drawn in 
Figure~\ref{K33}. A minor of a graph $G$ is a graph $H$ which can be 
obtained from $G$ by a sequence of vertex deletions, edge deletions and 
edge contractions. For more details about graph minors or the classical 
notations in graph theory used throughout the paper, the reader can refer to 
any general book on graph theory, for instance \cite{BM76,Die05}.

The main result of this paper is the following:

\begin{theorem}\label{main}
Let $G$ be a $3$-connected graph with no $K_{3,3}$ minor. Then every longest circuit of $G$ has a chord.
\end{theorem}

Denoting by $K^5$ the complete graph on five vertices, Kuratowski's 
theorem~\cite{Die05} 
states that a graph is planar if and only if it contains
neither $K^5$ nor $K_{3,3}$ as a minor. Therefore, Theorem~\ref{main} immediatly 
yields:

\begin{corollary}
Every longest circuit in a planar $3$-connected graph has a chord.
\end{corollary}

\medskip

\noindent {\bf Proof of Theorem~\ref{main} :}

Let $G$ be a  $3$-connected graph that does not contain $K_{3,3}$ as a minor 
and let $C$ be a longest circuit in $G$. Suppose that $C$ has no chord. We denote
by $v_0,\ldots ,v_p$ the vertices of $C$ in cyclic order.

Let $H_1,\ldots ,H_r$ be the connected components of $G\setminus C$. 
We denote by $N(i)$ the set of vertices of attachment of the component $H_i$, that is the
set of vertices of $C$ that are adjacent to a vertex of $H_i$.

Let $P$ be an arc of $C$, that is a connected subgraph of $C$, and let $i\in \{1,\ldots,r\}$. 
We say that $P$ is a {\em support} of $H_i$ if $N(i)\subset V(P)$.

\medskip

Let us first state some straightforward observations about the vertices of attachment 
of the components of $G\setminus C$:

\begin{lemma}\label{easy}
\begin{description}
\item[i)] For every $i\in \{1,\ldots,r\}$, $|N(i)|\geq 3$.
\item[ii)] For every $k\in \{0,\ldots,p\}$, there exist $i\in \{1,\ldots,r\}$ such that $v_k \in N(i)$.
\item[iii)] Two consecutive vertices on $C$ cannot belong to the same set $N(i)$.
\item[iv)] There are no integers $(k,l)\in \{0,\ldots,p\}^2$ and integers $(i,j)\in \{1,\ldots,r\}^2$ such that
  $\{v_k,v_l\} \subset N(i)$ and $\{v_{k+1},v_{l+1}\} \subset N(j)$. (By convention, $v_{p+1}=v_0$)
\end{description}
\end{lemma}

\pf

\begin{description}
\item [i)] This is a direct consequence of the $3$-connectedness of $G$.

\item[ii)] As $G$ is $3$-connected, $v_k$ has degree at least three. As $C$ has no chord, one of the neighbours
of $v_k$ therefore does not belong to $V(C)$.

\item[iii) and iv)] Figure~\ref{longest}
 show that both cases would contradict the fact that $C$ is a longest circuit of $G$.
Indeed, the bold circuits are longer, as their subpaths in the connected components of $G\setminus C$
contain at least one vertex.

\begin{figure} [htb]
\begin{center}
\psfrag{C}{$C$}
\psfrag{H1}{$H_i$}
\psfrag{H2}{$H_i$}
\psfrag{H3}{$H_j$}
\psfrag{v1}{$v_k$}
\psfrag{v2}{$v_{k+1}$}
\psfrag{v3}{$v_k$}
\psfrag{v4}{$v_{k+1}$}
\psfrag{v5}{$v_l$}
\psfrag{v6}{$v_{l+1}$}
\epsfig{file=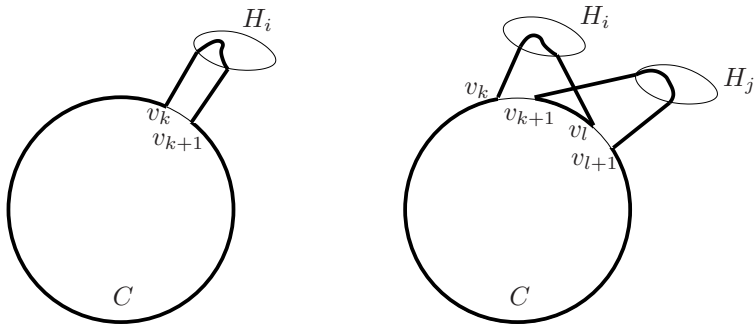, width=.8\textwidth}
\end{center}
\caption{Proof of the point iii) and iv) of Lemma~\ref{easy}}
\label{longest}
\end{figure}


\end{description}

\endpf

Consider an integer $i\in \{1,\ldots,r\}$ and a support $P$ of $H_i$ such that $P$ is minimal
with respect to inclusion among all the supports;
that is, there exists no pair $(j,Q)$ such that $Q$ is a support of $H_j$ and $Q$ is a proper subpath of $P$.
By reordering the vertices of $C$, we may assume that $P=v_0,\ldots ,v_q$, $2\leq q\leq p$.

We denote by $\Dot{P}$ the path $v_1\ldots v_{q-1}$. By Lemma~\ref{easy} ii) and iii), there exists at least one
index $j\neq i$
such that $N(j)\cap \Dot{P} \neq \emptyset$.

Let $j\neq i$ and $Q$ a support of $H_j$ be such that:
\begin{itemize}
\item $N(j)\cap \Dot{P} \neq \emptyset$
\item Subject to this condition, $P\cup Q$ is minimal (that is, there  exists no
 pair $(k,R)$ satisfying it and such that $P\cup R$ is a proper subpath 
of $P\cup Q$).
\item Subject to the two former conditions, $Q$ is minimal (that is, there  exists no
 pair $(k,R)$ with $N(k)\cap \Dot{P} \neq \emptyset$,  $P\cup R=P\cup Q$  and such that $R$ is a proper subpath 
of $Q$).
\end{itemize}

We can then choose six vertices, which will be denoted as the {\em vertices of interest}:
\begin{itemize}
\item Three distinct vertices of $N(i)$ in $P$, including both ends of $P$.
\item Three distinct vertices of $N(j)$ in $Q$, including both ends of $Q$ and  such that one of them is on $\Dot{P}$.
\end{itemize}

\begin{figure} [htb]
\begin{center}
\psfrag{i}{$i$}
\psfrag{j}{$j$}
\psfrag{k}{$k$}
\psfrag{P}{$P$}
\psfrag{Q}{$Q$}
\psfrag{Hk}{$H_k$}
\epsfig{file=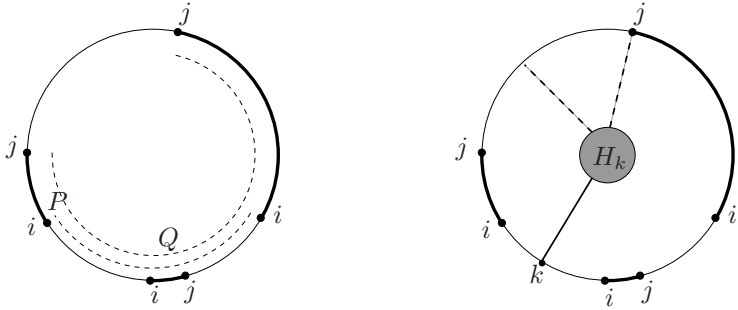, width=.8\textwidth}
\end{center}
\caption{The configuration when $P$ is a subpath of $Q$}
\label{PQ}
\end{figure}

Let us first assume that $P$ is a subpath of $Q$.
This case leads to the configuration shown in  Figure~\ref{PQ}. This figure has to been read in the following way:
\begin{itemize} 
\item The left hand side figure shows the repartition of the vertices of interest along $C$, where the bold arcs are possibly just
single vertices. To make the figure clearer, each 
vertex of interest is indicated by a label $i$ or $j$, meaning that it
belongs to $N(i)$ or $N(j)$, respectively.

\item The right hand side figure shows that in that configuration, $G$ admits a $K_{3,3}$ minor.

Indeed, 
 there exists an index $k\in \{1,\ldots,r\}$ different 
from $i$ and $j$ such that some vertex of $N(k)$ lies 
on a sub-arc of $P$ which is not bold: it is straightforward to verify that the absence of such a vertex would 
contradict part iii) or part iv) of Lemma~\ref{easy}.  

We can choose that vertex on the sub-arc shown on Figure~\ref{PQ} because if not,
there has to be a vertex labelled $j$ on that sub-arc and we obtain the same figure by symmetry.

Moreover, by minimality of $P\cup Q$, $N(k)$ either contains both ends of $P\cup Q$ or 
intersects the arc of $C$ with vertex set
$V(C)\setminus V(P\cup Q)$. In any case, one of the two dashed lines exists.
\medskip

Let then $A_1$, $A_2$ and $A_3$ be the three bold arcs of $C$. Let $B_1=H_i$, $B_2=H_j$ and $B_3$ be the union of $H_k$ and the arcs of $C\setminus (A_1\cup A_2\cup A_3)$ intersecting $N(k)$.
The above six sets are disjoint and connected, and each $A_i$ is linked to each $B_j$ by an edge.
Therefore, contracting each set to a single vertex, we obtain a $K_{3,3}$ minor of $G$.  

\end{itemize}

We use the same technique to deal with the remaining cases, that is when $P$ is not a subpath of $Q$. 

\begin{figure} [htb]
\begin{center}
\psfrag{(a)}{(a)}
\psfrag{(b)}{(b)}
\psfrag{(c)}{(c)}
\psfrag{i}{$i$}
\psfrag{j}{$j$}
\psfrag{k}{$k$}
\psfrag{P}{$P$}
\psfrag{Q}{$Q$}
\psfrag{Hk}{$H_k$}
\psfrag{or}{or}
\epsfig{file=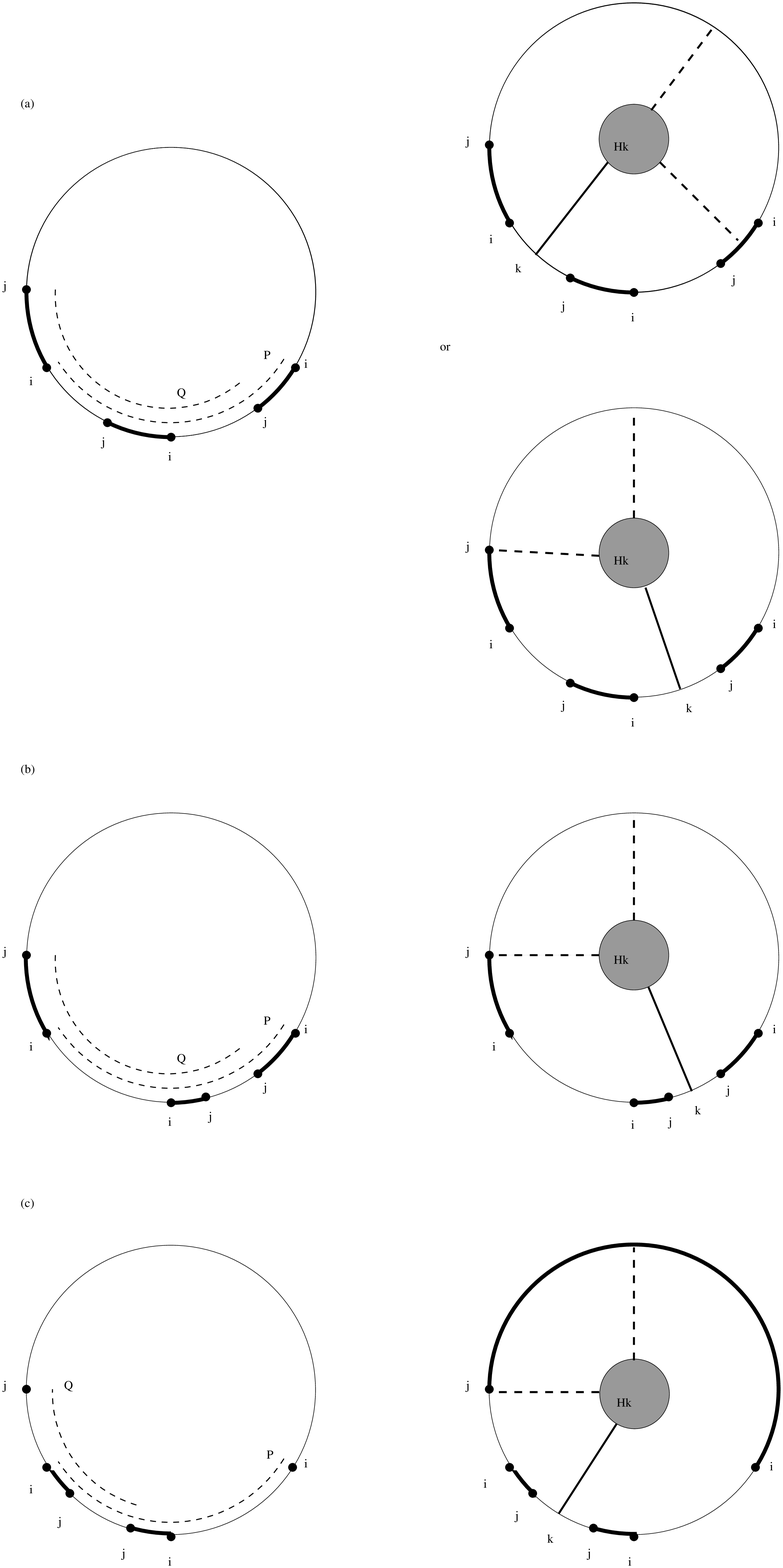, width=.8\textwidth}
\end{center}
\caption{The three possible configurations when $P$ in not a subpath of $Q$}
\label{PQ2}
\end{figure}

The possible orderings (up to symmetry) of the six vertices of interest are shown on the left side of Figure~\ref{PQ2}.
The right side has to be read as in Figure~\ref{PQ} and shows that in each case, one can exhibit a $K_{3,3}$ minor.

The three following remarks are to fully understand that Figure:
\begin{itemize}
\item To obtain the second dashed line on the top right figure, we use the minimality of $Q$ instead of
the minimality of $P\cup Q$  
\item In case (b), we can restrict ourselves to the case where the vertex labelled $k$ is between 
the two vertices of interest labelled $j$. Indeed, if not, there has to be a vertex labelled $i$ on that place
and we can use case (a) to conclude.
\item In case (c), the leftmost vertex labelled $j$ is distinct from the 
leftmost vertex labelled $i$ by minimality of $P$.

\end{itemize}

The fact that all possible configurations lead to a  contradiction ensures that $C$ has a chord.

\endpf

\bibliographystyle{plain}
\bibliography{Biblio}

\begin{thebibliography}{1}

\bibitem{Bon95}
J.A. Bondy.
\newblock Basic graph theory.
\newblock In L.~Lov\'asz M.~Gr\"{o}tschel and R.L. Graham, editors, {\em
  Handbook of Combinatorics}, pages 3--110. North-Holland, Amsterdam, 1995.

\bibitem{BM76}
J.A. Bondy and U.S.R. Murty.
\newblock {\em Graph Theory with Applications}.
\newblock American Elsevier, New-York, 1976.

\bibitem{Die05}
R.~Diestel.
\newblock {\em Graph Theory}.
\newblock Springer Verlag, Heidelberg, third edition edition, 2005.

\bibitem{KNZ07}
K.-I. Kawarabayashi, J.~Niu, and C.-Q. Zhang.
\newblock Chords of longest circuits in locally planar graphs.
\newblock {\em European J. Combin.}, 28(1):315--321, 2007.

\bibitem{LZ03}
X.~Li and C.-Q. Zhang.
\newblock Chords of longest circuits in $3$-connected graphs.
\newblock {\em Discrete Math.}, 268(1-3):199--206, 2003.

\bibitem{LZ03bis}
X.~Li and C.-Q. Zhang.
\newblock Chords of longest circuits of graphs embedded in torus and klein
  bottle.
\newblock {\em J. Graph Theory}, 43:1--23, 2003.

\bibitem{Tho84}
C.~Thomassen.
\newblock Plane representation of graphs.
\newblock In {\em Progress in Graph Theory (Waterloo, Ont., 1982)}, pages
  43--69. Academic Press, Toronto, ON, 1984.

\bibitem{Tho97}
C.~Thomassen.
\newblock Chords of longest cycles in cubic graphs.
\newblock {\em J. Combin. Theory Ser. B}, 71(2):211--214, 1997.

\bibitem{Zha87}
C.-Q. Zhang.
\newblock Longest cycles and their chords.
\newblock {\em J. Graph Theory}, 11(4):521--529, 1987.

\end{thebibliography}

\end{document}